\documentclass[12pt,reqno]{article}
\usepackage{amssymb,amscd,amsmath,amsthm,color}
\newtheorem{theorem}{Theorem}[section]

\newtheorem{cor}[theorem]{Corollary}

\usepackage{enumerate,amssymb}
\theoremstyle{definition}

\theoremstyle{remark}

\numberwithin{equation}{section}

\def\bK{\mathbb{K}}
\def\bM{\mathbb{M}}

\def\bB{\mathbb{B}}

\begin{document}
\baselineskip=15pt

\title{ Clarkson-McCarthy inequalities  with unitary and isometry orbits}

\author{ Jean-Christophe Bourin{\footnote{This research was supported by the French Investissements
 d'Avenir program, project ISITE-BFC (contract ANR-15-IDEX-03).
}} \ and  Eun-Young Lee{\footnote{This research was supported by
Basic Science Research Program through the National Research
Foundation of Korea (NRF) funded by the Ministry of
Education (NRF-2018R1D1A3B07043682)}  }}

\date{ }

\maketitle

\vskip 10pt\noindent
{\small
{\bf Abstract.}  A refinement of a  trace inequality of McCarthy establishing the uniform convexity of the Schatten $p$-classes for $p>2$ is proved: if $A,B$ are two $n$-by-$n$ matrices, then there exists some pair of  $n$-by-$n$ unitary matrices $U,V$ such that
$$
U\left| \frac{A+B}{2}\right|^pU^* +V\left| \frac{A-B}{2}\right|^pV^*\le \frac{|A|^p+|B|^p}{2}.
$$
A similar statement holds for compact Hilbert space operators. Another improvement of McCarthy's inequality is given via the new operator parallelogramm law,
\begin{equation*}
\left|A+B\right|^2\oplus\left|A-B\right|^2= U_0(|A|^2+|B|^2)U_0^*+V_0(|A|^2+|B|^2)V_0^*
\end{equation*}
for some pair of $2n$-by-$n$ isometry matrices $U_0,V_0$.

\vskip 5pt\noindent
{\it Keywords.} Matrix inequalities,    Unitary orbits, Clarkson-McCarthy inequalities.
\vskip 5pt\noindent
{\it 2010 mathematics subject classification.} 47A30, 15A60.
}

\section{Introduction } 

Let $\bM_n$ denote the space of complex $n\times n$ matrices and let $\bM_n^+$ stand for the positive (semi-definite) cone. In \cite{McC}  McCarthy established a  Clarkson type inequality for the Schatten $p$-norms $\|\cdot\|_p$ and $A,B\in\bM_n$,
\begin{equation}\label{eq1}
2\left(\| A\|_p^p +\| B\|_p^p\right) \le \| A+B\|_p^p +  \| A-B\|_p^p
\le 
2^{p-1}\left(\| A\|_p^p +\| B\|_p^p\right), \quad p\ge 2.
\end{equation}
For $0< p\le 2$, these inequalities are reversed. In the commutative (or scalar) case, these inequalities are almost trivial.
Note that the left hand side inequality is equivalent to the right hand side one by $A=X+Y$, $B=X-Y$. Those inequalities for $p\ge 2$ show that the unit ball  for the Schatten $p$-norm is uniformly convex, this is more readily apparent when written as
\begin{equation}\label{Clark}
 \left\| \frac{ A+B}{2}\right\|_p^p +   \left\| \frac{ A-B}{2}\right\|_p^p
\le 
\frac{\| A\|_p^p +\| B\|_p^p}{2}, \quad p\ge 2,
\end{equation}
i.e.,
\begin{equation}\label{Clark'}
{\mathrm{Tr\,}} \left| \frac{ A+B}{2}\right|^p +   {\mathrm{Tr\,}} \left| \frac{ A-B}{2}\right|^p
\le 
\frac{{\mathrm{Tr\,}} | A|^p +{\mathrm{Tr\,}} | B|^p}{2}, \quad p\ge 2.
\end{equation}
Thus, if $\|A\|_p=\|B\|_p=1$ and $\| A-B\|_p = \varepsilon$, 
$$
\left\| \frac{ A+B}{2}\right\|_p \le \left(1-(\varepsilon/2)^p\right)^{1/p}, \quad  p\ge 2,
$$
which estimates the uniform convexity modulus of the Schatten $p$-classes for $p\ge 2$. There also exists a Clarkson type inequality showing the uniform convexity of the Schatten $p$-classes in case of $1<p<2$. This case is not as simple as the case $p>2$ and a Three Lines Theorem argument is required.  It seems that no real analytic proof are known (the original proof given by McCarthy collapses, see \cite{FK}, p.\ 297)

The Clarkson-McCarthy inequalties \eqref{eq1}  have been nicely extended to a large class of unitarily invariant norms by Bhatia and Holbrook \cite{BH}. Recall that a unitarily invariant norm on $\bM_n$, also called a symmetric norm, satisfies $\|UAV\|=\|A\|$ for all 
$A\in\bM_n$ and all unitary matrices $U,V\in\bM_n$. Another remarkable generalization  due to Hirzallah and Kittaneh \cite{HK} states that
$$
2\left\| |A|^p+|B|^p\right\| \le \left\| |A+B|^p +  |A-B|^p \right\|
\le
2^{p-1}\left\| |A|^p +|B|^p \right\|, \quad p\ge 2,
$$
for all symmetric norms, and these inequalities are reversed for $0<p\le 2$. This double inequality can be rewritten in a similar form as \eqref{Clark} and, thanks to Fan's dominance principle, this can be subsumed as the weak majorization
\begin{equation}\label{Clark2}
 \sum_{j=1}^k\lambda_{j}^{\downarrow}\left(\left| \frac{ A+B}{2}\right|^p +   \left| \frac{ A-B}{2}\right|^p\right)
\le 
 \sum_{j=1}^k\lambda_{j}^{\downarrow}\left(\frac{| A|^p +| B|^p}{2}\right), \quad k=1,2,\ldots
\end{equation}
where $\lambda_{1}^{\downarrow}(X )\ge \lambda_{2}^{\downarrow}(X )\ge \cdots\ge \lambda_{n}^{\downarrow}(X ) $ stand for the eigenvalues of $X\in\bM_n^+$. For $k=n$, we recapture \eqref{Clark'}, hence \eqref{Clark2} is a considerable improvement of \eqref{Clark}.

We will show in the next two sections some refinements of \eqref{Clark}. These results yields some eigenvalue estimates which complete \eqref{Clark2}. The last section show how to modify the statements of Section 2 when dealing with operators
 on an infinite dimensional Hilbert space.

\section{Unitary orbits}

The operator parallelogramm law does hold, for every pair $A,B\in\bM_n$,
\begin{equation}\label{parall}
|A+B|^2+ |A-B|^2 =2\left(|A|^2+|B|^2\right),
\end{equation}
equivalently,
$$
\left| \frac{A+B}{2}\right|^2 +\left| \frac{A-B}{2}\right|^2 = \frac{|A|^2+|B|^2}{2}.
$$
 We will state two theorems which show for $p\neq 2$ how this equality is modified as operator inequalities involving unitary orbits. Taking traces we recapture \eqref{eq1}-\eqref{Clark}.

\vskip 5pt
\begin{theorem}\label{th1} Let $A,B\in\bM_n$ and $p>2$. Then there exist two unitaries $U,V\in\bM_n$ such that
$$
U\left| \frac{A+B}{2}\right|^pU^* +V\left| \frac{A-B}{2}\right|^pV^*\le \frac{|A|^p+|B|^p}{2}.
$$
\end{theorem}

\vskip 5pt
\begin{proof} Note that
$$
\left| \frac{A+B}{2}\right|^p= \left(\frac{|A|^2+|B|^2 +A^*B +B^*A}{4}\right)^{p/2}
$$
and 
$$
\left| \frac{A-B}{2}\right|^p= \left(\frac{|A|^2+|B|^2 -(A^*B +B^*A)}{4}\right)^{p/2}.
$$
Now, recall  \cite[Corollary 3.2]{BL}: Given two positive matrices $X,Y$ and a monotone convex function $g(t$) defined on $[0,\infty)$ such that $g(0)\le 0$, we have 
\begin{equation}\label{key1}
g(X+Y) \ge U_0g(X)U_0* + V_0g(Y)V_0*
\end{equation}
for some pair of unitary matrices $U_0$ and $V_0$. Applying this to $g(t)=t^{p/2}$,
$$
X=\frac{|A|^2+|B|^2 +A^*B +B^*A}{4}
$$
and
$$
Y=\frac{|A|^2+|B|^2 -(A^*B +B^*A)}{4},
$$
 we obtain
\begin{equation}\label{e1}
\left(\frac{|A|^2+|B|^2}{2}\right)^{p/2} \ge U_0\left| \frac{A+B}{2}\right|^pU_0^* +V_0\left| \frac{A-B}{2}\right|^pV_0^*.
\end{equation}
Next, recall \cite[Corollary 2.2]{BL}:  Given two positive matrices $X,Y$ and a monotone convex function $g(t$) defined on $[0,\infty)$ , we have 
\begin{equation}\label{key2}
\frac{g(X)+g(Y)}{2} \ge Wg\left(\frac{X+Y}{2}\right)W^*
\end{equation}
for some unitary matrix $W$. Applying this to $g(t)=t^{p/2}$, $X=|A|^2$ and $Y=|B|^2$, we get
\begin{equation}\label{e2}
 \frac{|A|^p+|B|^p}{2}\ge W\left(\frac{|A|^2+|B|^2}{2}\right)^{p/2}W^*.
\end{equation}
Combining \eqref{e1} and \eqref{e2} completes the proof with $U=WU_0$ and $V=WV_0$.
\end{proof}

\vskip 5pt
\begin{cor}\label{cor1}
 Let $A,B\in\bM_n$ and $p>2$. Then, for all $k=1,2,\ldots,n$,
$$
\sum_{j=1}^k \lambda_j^{\uparrow}\left(\frac{|A|^p+|B|^p}{2}\right) \ge 
\sum_{j=1}^k \lambda_j^{\uparrow}\left(\left|\frac{A+B}{2}\right|^p\right)
+
\sum_{j=1}^k \lambda_j^{\uparrow}\left(\left|\frac{A-B}{2}\right|^p\right).
$$
\end{cor}

\vskip 5pt
\begin{cor}\label{cor2}
 Let $A,B\in\bM_n$ and $p>2$. Then, for all $k=1,2,\ldots,n$,
$$
\left\{\prod_{j=1}^k \lambda_j^{\uparrow}\left(\frac{|A|^p+|B|^p}{2}\right)\right\}^{1/k} \ge 
\left\{\prod_{j=1}^k \lambda_j^{\uparrow}\left(\left|\frac{A+B}{2}\right|^p\right)\right\}^{1/k}
+
\left\{\prod_{j=1}^k \lambda_j^{\uparrow}\left(\left|\frac{A-B}{2}\right|^p\right)\right\}^{1/k}.
$$
\end{cor}

\vskip 5pt
Here  $\lambda_{1}^{\uparrow}(X )\le \lambda_{2}^{\uparrow}(X )\le \cdots\le \lambda_{n}^{\uparrow}(X ) $ stand for the eigenvalues of $X\in\bM_n^+$ arranged in the nondecresaing order.
These two corollaries follow from the theorem and the fact that the functionals on $\bM_n^+$
$$X\mapsto \sum_{j=1}^k \lambda_j^{\uparrow}(X)$$
and
$$X\mapsto\left\{\prod_{j=1}^k  \lambda_j^{\uparrow}(X)\right\}^{1/k}$$
are two basic examples of symmetric anti-norms, see \cite{BouHiai1}.

The next corollary follows from Theorem \ref{th1} combined with a classical inequality of Weyl \cite[page 53]{Z} for the eigenvalues of the sum of two Hermitian matrices.

\vskip 5pt
\begin{cor}\label{cor3}
 Let $A,B\in\bM_n$ and $p>2$. Then, for all $j,k\in\{0,\ldots,n-1\}$ such that $j+k+1\le n$,
$$
 \lambda_{j+1}^{\downarrow}\left(\frac{|A|^p+|B|^p}{2}\right) \ge 
 \lambda_{j+k+1}^{\downarrow}\left(\left|\frac{A+B}{2}\right|^p\right)
+
 \lambda_{k+1}^{\uparrow}\left(\left|\frac{A-B}{2}\right|^p\right).
$$
\end{cor}

\vskip 5pt
For a monotone concave function $g(t$) defined on $[0,\infty)$ such that $g(0)\ge 0$, the inequalities \eqref{key1} and \eqref{key2} are reversed.
Applying this to $g(t)=t^{q/2}$, $2>q>0$, the same proof than that of Theorem \ref{th1} gives the following statement.

\vskip 5pt
\begin{theorem}\label{th2} Let $A,B\in\bM_n$ and $2>q>0$. Then, for some unitaries $U,V\in\bM_n$,
$$
U\left| \frac{A+B}{2}\right|^qU^* +V\left| \frac{A-B}{2}\right|^qV^*\ge \frac{|A|^q+|B|^q}{2}.
$$
\end{theorem}

\vskip 5pt
By using Weyl's inequality, this theorem yields an interesting eigenvalue estimate.

\vskip 5pt
\begin{cor}\label{cor4}
 Let $A,B\in\bM_n$ and $2>q>0$. Then, for all $j,k\in\{0,\ldots,n-1\}$ such that $j+k+1\le n$,
$$
 \lambda_{j+k+1}^{\downarrow}\left(\frac{|A|^q+|B|^q}{2}\right) \le 
 \lambda_{j+1}^{\downarrow}\left(\left|\frac{A+B}{2}\right|^q\right)
+
 \lambda_{k+1}^{\downarrow}\left(\left|\frac{A-B}{2}\right|^q\right).
$$
\end{cor}

\section{Another parallelogram law}

We shall point out a new operator parallelogram law where the usual sum in \eqref{parall} is replaced by a direct sum.
 A matrix $V\in \bM_{m,n}$, the space of $m\times n$ matrices, is called an {\it isometry} whenever $V^*V$ is the identity of $\bM_n$. 

\vskip 5pt
\begin{theorem} \label{th5} Let $A,B\in\bM_n$. Then, for some isometries  $U,V\in\bM_{2n,n}$,
\begin{equation*}
\left|A+B\right|^2\oplus\left|A-B\right|^2= U(|A|^2+|B|^2)U^*+V(|A|^2+|B|^2)V^*
\end{equation*}
\end{theorem}

\vskip 5pt
\begin{proof} Note that
$$
\begin{pmatrix}
A&B \\ B&A
\end{pmatrix}
$$
is unitarily equivalent to 
$$
\begin{pmatrix}
A+B&0 \\ 0&A-B
\end{pmatrix}
$$
via the unitary congruence implemented by
$$
W={\frac{1}{\sqrt{2}}}
\begin{pmatrix}
I&I \\ I&-I
\end{pmatrix}.
$$
Thus
$\left|A+B\right|^2\oplus\left|A-B\right|^2$
is unitarily equivalent to
\begin{equation}\label{block}
\left|\begin{pmatrix}
A&B \\ B&A
\end{pmatrix}\right|^2=\begin{pmatrix}
|A|^2+|B|^2&A^*B +B^*A \\ A^*B +B^*A &|A|^2+|B|^2
\end{pmatrix}.
\end{equation}
Now, recall \cite[Lemma 3.4]{BL}: Given any positive matrix partitionned in four $n$-by-$n$ blocks, we can find two unitary matrices $U,V\in\bM_{2n}$ such that
$$
\begin{pmatrix}
X&Y \\ Y^*&Z
\end{pmatrix}=U\begin{pmatrix}
X&0 \\ 0 &0
\end{pmatrix}U^*+
V_0\begin{pmatrix}
0&0 \\ 0 &Z
\end{pmatrix}V^*.
$$
Applying this to \eqref{block}
 we then have two unitary matrices $U,V\in\bM_{2n}$ such that 
$$\left|A+B\right|^2\oplus\left|A-B\right|^2=U\begin{pmatrix}|A|^2+|B|^2&0\\ 0&0
\end{pmatrix}U^* + 
V\begin{pmatrix}0&0\\ 0&|A|^2+|B|^2
\end{pmatrix}V^*
$$
which is a statement equivalent to our theorem.
\end{proof}

\vskip 5pt
The next corollary is a matrix version of the  identity $|z|^2=x^2+y^2$ for complex numbers $z=x+iy$. It follows from  Theorem \ref{th5} applied  to $A=X$ and $B=iY$.

\vskip 5pt
\begin{cor} \label{corabs} Let $Z\in\bM_n$ with Cartesian decomposition $Z=X+iY$. Then, for some isometries  $U,V\in\bM_{2n,n}$,
\begin{equation*}
\left|Z\right|^2\oplus\left|Z\right|^2= U(X^2+Y^2)U^*+V(X^2+Y^2)V^*.
\end{equation*}
\end{cor}

\vskip 5pt
From this and the reverse form of \eqref{key1} for $g(t)=\sqrt{t}$ we see that the scalar identity $|z|=\sqrt{x^2+y^2}$ is transformed 
 in the matrix setting as an operator inequality.

\vskip 5pt
\begin{cor} \label{corabs2} Let $Z\in\bM_n$ with Cartesian decomposition $Z=X+iY$. Then, for some isometries  $U,V\in\bM_{2n,n}$,
\begin{equation*}
\left|Z\right|\oplus\left|Z\right| \le  U\sqrt{X^2+Y^2}\,U^*+V\sqrt{X^2+Y^2}\,V^*.
\end{equation*}
\end{cor}

\vskip 5pt
By applying twice the  inequality \eqref{key1} to the convex function $g(t)=t^{p/2}$, our parallelogramm law
yields the next corollary.

\vskip 5pt
\begin{cor} \label{cordirect1} Let $A,B\in\bM_n$ and $p>2$. Then, for some isometries  $U_0,V_0,U_1,V_1\in\bM_{2n,n}$,
\begin{equation*}
\left|A+B\right|^p\oplus\left|A-B\right|^p\ge U_0|A|^pU_0^*+V_0|B|^pV_0^*+ U_1|A|^pU_1^*+V_1|B|^pV_1^*.
\end{equation*}
\end{cor}

\vskip 5pt
If one uses the concave function $g(t)=t^{q/2}$, we obtain a reversed inequality.

\vskip 5pt
\begin{cor} \label{cordirect2} Let $A,B\in\bM_n$ and $2>q>0$. Then, for some isometries  $U_0,V_0,U_1,V_1\in\bM_{2n,n}$,
\begin{equation*}
\left|A+B\right|^q\oplus\left|A-B\right|^q\le  U_0|A|^qU_0^*+V_0|B|^qV_0^*+ U_1|A|^qU_1^*+V_1|B|^qV_1^*.
\end{equation*}
\end{cor}

\vskip 5pt
Our parallelogram law yields another extension of the Clarkson-McCarthy trace inequality \eqref{Clark}.

\vskip 5pt
\begin{theorem}\label{thdirect2} Let $A,B\in\bM_n$ and $p>2$. Then, for some isometries  $U,V\in\bM_{2n,n}$,
$$
\left| \frac{A+B}{2}\right|^p \oplus \left| \frac{A-B}{2}\right|^p\le \frac{1}{2}\left\{U\frac{|A|^p+|B|^p}{2}U^* +V\frac{|A|^p+|B|^p}{2}V^* \right\}.
$$
\end{theorem}

\vskip 5pt
\begin{proof} Theorem \ref{th5} says that
\begin{equation*}
\left|\frac{A+B}{2}\right|^2\oplus\left|\frac{A-B}{2}\right|^2= \frac{1}{2}\left\{U\frac{|A|^2+|B|^2}{2}U^*+V\frac{|A|^2+|B|^2}{2}V^*\right\}
\end{equation*}
and applying twice inequality \eqref{key2} with the convex function $g(t)=t^{p/2}$ completes the proof.
\end{proof}

\vskip 5pt
The simplest form of Weyl's inequality,  $\lambda_{2j+1}^{\downarrow}(S+T) \le  \lambda_{j+1}^{\downarrow}(S)+ \lambda_{j+1}^{\downarrow}(T)$ for $S,T\in\bM_n^+$, and
Theorem \ref{thdirect2} provide one more Clarkson-McCarthy type eigenvalue estimate.

\vskip 5pt
\begin{cor}\label{cor5}
 Let $A,B\in\bM_n$ and $p>2$. Then, for all $j,\in\{0,\ldots,n-1\}$,
$$
 \lambda_{2j+1}^{\downarrow}\left(\left| \frac{A+B}{2}\right|^p \oplus \left| \frac{A-B}{2}\right|^p\right) \le 
 \lambda_{j+1}^{\downarrow}\left(  \frac{|A|^p+|B|^p}{2} \right).
$$
\end{cor}

\vskip 5pt
The reverse form of Theorem \ref{thdirect2} occurs for $2>q>0$.

\vskip 5pt
\begin{theorem}\label{thdirect3} Let $A,B\in\bM_n$ and $2>q>0$. Then, for some isometries  $U,V\in\bM_{2n,n}$,
$$
\left| \frac{A+B}{2}\right|^q \oplus \left| \frac{A-B}{2}\right|^q\ge \frac{1}{2}\left\{U\frac{|A|^q+|B|^q}{2}U^* +V\frac{|A|^q+|B|^q}{2}V^* \right\}.
$$
\end{theorem}

\section{Infinite dimension}

The Schatten classes are usually defined as classes of compact operators on an infinite dimensional, separable Hilbert space.
In this setting, a slight modification of the previous theorems is necessary. The correct statements for compact operators require to use partial isometries rather than unitary operators.

Let $g(t)$ be an increasing continuous function defined on $[0,\infty)$ such that $g(0)$. Then the functional calculus  $g:\bK^+\to \bK^+$, $X\mapsto g(X)$
is norm continuous. Thus, if $X_n$, $n\ge 1$, is a sequence of finite rank operators in $\bK^+$ converging to $X$, then we have
$$\lambda_j^{\downarrow}(X_n)\to \lambda_j^{\downarrow}(X), \quad j=1,2,\ldots$$
and
$$\lambda_j^{\downarrow}(g(X_n))\to \lambda_j^{\downarrow}(g(X)), \quad j=1,2,\ldots$$
as $n\to \infty$. Now let $X,Y\in\bK^+$ and pick two sequences of finite rank positive operators $X_n$ and $Y_n$ such that $X_n\to X$
and $Y_n\to Y$. Since we deal with finite rank operators, if $g(t)$ is convex,  \eqref{key2} yields
\begin{equation*}
\frac{g(X_n)+g(Y_n)}{2} \ge Wg\left(\frac{X_n+Y_n}{2}\right)W^*
\end{equation*}
for some unitary operator $W$. This says that
\begin{equation*}
\lambda_j^{\downarrow}\left(\frac{g(X_n)+g(Y_n)}{2}\right) \ge \lambda_j^{\downarrow}\left(g\left(\frac{X_n+Y_n}{2}\right)\right), \quad j=1,2,\ldots,
\end{equation*}
 and, letting $n\to\infty$, we get
\begin{equation*}
\lambda_j^{\downarrow}\left(\frac{g(X)+g(Y)}{2}\right) \ge \lambda_j^{\downarrow}\left(g\left(\frac{X+Y}{2}\right)\right), \quad j=1,2,\ldots,
\end{equation*}
which exactly says that
\begin{equation*}
\frac{g(X)+g(Y)}{2} \ge Wg\left(\frac{X+Y}{2}\right)W^*
\end{equation*}
for some partial isometry $W$ such that ${\mathrm{supp}}(W)={\mathrm{supp}}(X+Y)$.

Therefore, for compact operators, Equation \eqref{key2} of the proof of Theorem \ref{th1} still holds with a partial isometry  $W$ such that ${\mathrm{supp}}(W)={\mathrm{supp}}(X+Y)$.
A similar statement holds for Equation  \eqref{key1}: Given $X,Y\in\bK^+$ and an increasing, strictly convex function $g(t)$ such that $g(0)$, it is shown in \cite[Theorem 2.1(i)]{BHL} that
\begin{equation*}
g(X+Y) \ge U_0g(X)U_0^* + V_0g(Y)V_0^*
\end{equation*}
for some partial isometries $U_0,V_0$ such that ${\mathrm{supp}}(U_0)={\mathrm{supp}}(X)$ and
${\mathrm{supp}}(V_0)={\mathrm{supp}}(Y)$. Hence, we may brought the proof of Theorem \ref{th1} to the setting of compact operators and we obtain the following refinement of the Clarkson-McCarthy inequalities for Schatten $p$-classes, $p>2$.

\vskip 5pt
\begin{theorem}\label{th3} Let $A,B\in\bK$ and $p>2$. Then there exist two partial isometries $U,V\in\bB$ such that ${\mathrm{supp}}(U)= {\mathrm{supp}}(A+B)$,  ${\mathrm{supp}}(V)= {\mathrm{supp}}(A-B)$, and
$$
U\left| \frac{A+B}{2}\right|^pU^* +V\left| \frac{A-B}{2}\right|^pV^*\le \frac{|A|^p+|B|^p}{2}.
$$
\end{theorem}

\vskip 5pt
We may also give some versions of Theorems \ref{th1} and \ref{th2} for measurable operators affiliated to a semi-finite  von Neumann algebra (we refer to \cite{FK} for a background on measurable operators). We close the paper by stating such a result for a type ${\mathrm{II}}_1$ factor. We do not give details, we just mention that the result  follows 
from \cite[Lemma 3.3]{BouHiai2} and \cite[Theorem 3.2]{BouHiai2}. 

We denote by $\overline{\mathcal{N}}$ the space of measurable operators affiliated to a type  ${\mathrm{II}}_1$  factor ${\mathcal{N}}$.

\vskip 5pt
\begin{theorem}\label{th4}
Let  $A,B\in\overline{\mathcal{N}}$ and $2>q>0$. Then, for every $\varepsilon>0$, there exist two unitaries $U,V\in{\mathcal{N}}$
such that
$$
U\left| \frac{A+B}{2}\right|^qU^* +V\left| \frac{A-B}{2}\right|^qV^*\ge \frac{|A|^q+|B|^q}{2}-\varepsilon I.
$$
\end{theorem}

\vskip 5pt
For $q>2$, the reverse inequality holds, with an $\varepsilon I$ term instead of  $-\varepsilon I$.

\vskip 5pt
\noindent
Laboratoire de math\'ematiques, 

\noindent
Universit\'e de Bourgogne Franche-Comt\'e, 

\noindent
25 000 Besan\c{c}on, France.

\noindent
Email: jcbourin@univ-fcomte.fr

  \vskip 10pt
\noindent Department of mathematics, KNU-Center for Nonlinear
Dynamics,

\noindent
Kyungpook National University,

\noindent
 Daegu 702-701, Korea.

\noindent Email: eylee89@knu.ac.kr

\end{document}